\documentclass{article}%
\usepackage{amsfonts}
\usepackage{amsmath}
\usepackage{amssymb}
\usepackage{graphicx}%
\newcommand{\bpartial}{\mathop{\partial\kern -4pt\raisebox{.8pt}{$|$}}}
\newcommand{\bra}{\mathopen{[\kern-1.6pt[}}
\newcommand{\ket}{\mathclose{]\kern-1.5pt]}}
\newcommand{\bbra}{\mathopen{[\kern-2.2pt[\kern-2.3pt[}}
\newcommand{\bket}{\mathclose{]\kern-2.1pt]\kern-2.3pt]}}
\newcommand{\slg}{\mbox{\bfseries\slshape g}}

\makeindex
\begin{document}

\title{Multivector and Extensor Fields on Smooth Manifolds}
\author{{\footnotesize A. M. Moya}$^{2},${\footnotesize V. V. Fern\'{a}ndez}$^{1}%
${\footnotesize , and W. A. Rodrigues Jr}.$^{1}${\footnotesize \ }\\$^{1}\hspace{-0.1cm}${\footnotesize Institute of Mathematics, Statistics and
Scientific Computation}\\{\footnotesize \ IMECC-UNICAMP CP 6065}\\{\footnotesize \ 13083-859 Campinas, SP, Brazil }\\{\footnotesize e-mail: walrod@ime.unicamp.br virginvelfe@accessplus.com.ar}\\{\footnotesize \ }$^{2}${\footnotesize Department of Mathematics, University
of Antofagasta, Antofagasta, Chile} \\{\footnotesize e-mail: mmoya@uantof.cl} }
\maketitle

\begin{abstract}
The main objective of this paper (second in a series of four) is to show how
the Clifford and extensor algebras methods introduced in a previous paper of
the series are indeed powerful tools for performing sophisticated calculations
appearing in the study of the differential geometry of a $n$-dimensional
manifold $M$ of arbitrary topology, supporting a metric field $%
%TCIMACRO{\TeXButton{slg}{\slg}}%
%BeginExpansion
\slg
%EndExpansion
$ (of given signature $(p,q)$) and an arbitrary connection $\nabla$.
Specifically, we deal here with the theory of multivector and extensor fields
on $M$. Our approach does not suffer the problems of earlier attempts which
are restricted to vector\emph{\ }manifolds. It is based on the existence of
canonical algebraic structures over the canonical (vector) space associated to
a local chart $(U_{o},\phi_{o})$ of a given atlas of $M$. The key concepts of
$a$-directional ordinary derivatives of multivector and extensor fields are
defined and their properties studied. Also, we recall the Lie algebra of
smooth vector fields in our formalism, the concept of Hestenes derivatives and
present some illustrative applications.

\end{abstract}
\tableofcontents

\section{Introduction}

The main purpose of the present paper, the second, in a series of four (and
sequel papers in the series) is to show how the Clifford and extensor algebras
methods developed in a previous paper \cite{3} are indeed a powerful tool for
doing calculations appearing in the study of the differential geometry of a
$n$-dimensional smooth manifold $M$ of arbitrary topology, equipped supporting
a metric field $%
%TCIMACRO{\TeXButton{slg}{\slg}}%
%BeginExpansion
\slg
%EndExpansion
$ \ (of signature $(p,q)$, $p+q=n$) and an arbitrary connection $\nabla$. Of
course, a sophisticated way to apply Clifford and extensor algebras methods to
the study of the differential geometry of manifolds is the Clifford bundle
formalism based on the Clifford bundle of multivector fields\footnote{Or,
which is more appropriate, the Clifford bundle of multiforms
\cite{rodoliv2006}.} $\mathcal{C\ell}(TM,g)$, which is described, e.g., in
\cite{rodoliv2006}. However, here we want to avoid (as much as possible) the
theory of vector and principal bundles and connections on them, and indeed our
intention is to give a simple approach to the subject that can be ready used
by interested physicists. So, to begin, we recall that any effective practical
calculation in differential geometry, starts with the selection of an
appropriate local chart, say $(U_{o},\phi_{o})$ with coordinates $\{x^{\mu}\}$
of a given atlas of $M$ . Instead of working directly with multivector and
extensor fields\footnote{Which are sections of appropriate vector bundles.
See, e.g., \cite{rodoliv2006}}, the main idea of our approach consists in
working with the \textit{representatives }of those \textit{fields }on\textit{
}$U_{o}\subset M$\textit{ }through the use of canonical algebraic structures
constructed over a canonical vector space\footnote{It thus does not suffer the
problems of some earlier attempts that have been restricted to vector
manifolds \cite{6,7}.} associate to a local chart $(U_{o},\phi_{o})$. These
key concepts are introduced in Section 2. In Section 3 we introduce the $a$
-directional \textit{ordinary} derivative $a\cdot\partial_{o}$\ of multivector
fields\footnote{Which is indeed a connection on $U_{o}\subset M$.}, and recall
the construction of the Lie algebra of smooth vector fields within our
formalism. We introduce also the concept of Hestenes derivatives of
multivector field $X$ , $\partial_{o}\ast X$, which will play an important
role in sequel papers of this series. The $a$-directional derivatives of
extensor fields is introduced in Section 4. Properties of those derivative
operators are proved in details. In Section 5 we show how to use (if
necessary) a different coordinate system associated to a local chart
$(U,\phi)$, for $U\cap U_{o}\subset U_{o}$, besides the one associate to the
local chart $(U_{o},\phi_{o})$ originally selected to define the canonical
space. This defines a new $a$-directional derivative, denoted $a\cdot\partial
$. We present also some elementary applications of the formalism, introducing
covariant and contravariant frame fields associated to a given coordinate
system and also the notion of Jacobian fields, which permit us to find the
relation between $a\cdot\partial_{o}$ and $a\cdot\partial$.

\section{Canonical Space}

Let $M$ be a $n$-dimensional smooth manifold. As well known, associated to any
point $o\in M$ there always exists a\emph{\ local chart} $(U_{o},\phi_{o})$ of
a given atlas of $M$ such that $o\in U_{o}$ and $\phi_{o}(o)=(0,\ldots,0)$.
Recall that there exist exactly $n$ scalar functions $\phi_{o}^{\mu}%
:U_{o}\rightarrow\mathbb{R}$ such that $\phi_{o}(p)=(\phi_{o}^{1}%
(p),\ldots,\phi_{o}^{n}(p))$ which are the coordinate functions of
$(U_{o},\phi_{o}).$ We have that $\phi_{o}^{\mu}=\pi^{\mu}\circ\phi_{o},$
where $\pi^{\mu}$ is the well-known $\mu$-\emph{projection mapping}%
\footnote{Recall that $\pi^{\mu}:\mathbb{R}^{n}\rightarrow\mathbb{R}$ is such
that if $a=(a^{1},\ldots,a^{n})\in\mathbb{R}^{n},$ then $\pi^{\mu}(a)=a^{\mu
}\in\mathbb{R}$ for each $\mu=1,\ldots,n.$ Note that $a=(\pi^{1}(a),\ldots
,\pi^{n}(a)).$} of $\mathbb{R}^{n}.$ Any \emph{point} $p\in U_{o}$ is then
localized by a $n$-uple of \emph{real numbers} $\phi_{o}(p)\in\mathbb{R}^{n}$,
and the $n$ real numbers
\begin{equation}
x_{o}^{\mu}=\phi_{o}^{\mu}(p),\text{ for each }\mu=1,\ldots,n \label{CVS.1}%
\end{equation}
are the \emph{position} \emph{coordinates} of $p$ with respect to $(U_{o}%
,\phi_{o}).$

As well known, the $n$\emph{\ coordinate} \emph{tangent vectors}\footnote{A
\emph{tangent vector} at $p\in U_{o}$ associated to the $\mu$-th position
coordinate $x_{o}^{\mu},$ namely $\left.  \dfrac{\partial}{\partial x_{o}%
^{\mu}}\right\vert _{p},$ can be defined by $\left.  \dfrac{\partial}{\partial
x_{o}^{\mu}}\right\vert _{p}f=(\dfrac{\partial}{\partial x_{o}^{\mu}}%
f\circ\phi_{o}^{-1})\circ\phi_{o}(p),$ for all $f\in\mathcal{C}^{\infty}(p).$%
}\emph{\ }$\left.  \dfrac{\partial}{\partial x_{o}^{1}}\right\vert _{p}%
,\ldots,\left.  \dfrac{\partial}{\partial x_{o}^{n}}\right\vert _{p}$ define a
\emph{natural basis} at each point $p\in U_{o}$ for the \emph{tangent space}
$T_{p}M$, called \emph{coordinate vector basis }at $p\in U_{o},$ i.e.,
\begin{equation}
\{\left.  \dfrac{\partial}{\partial x_{o}^{\mu}}\right\vert _{p}\}.
\label{CVS.2}%
\end{equation}
We recall, in order to fix our notations that for any $\mathbf{v}_{p}\in
T_{p}M$ we have the elementary expansion%
\begin{equation}
\mathbf{v}_{p}=\mathbf{v}_{p}\phi_{o}^{\mu}\left.  \dfrac{\partial}{\partial
x_{o}^{\mu}}\right\vert _{p},
\end{equation}
where $\mathbf{v}_{p}\phi_{o}^{\mu}$ denotes the components of $\mathbf{v}%
_{p}$ in the natural basis.

To continue, we introduce an equivalence relation on the set of the tangent
vectors on $\underset{p\in U_{o}}{\bigcup}T_{p}M$ as follows. Let
$\mathbf{v}_{a}\in T_{a}M$ and $\mathbf{v}_{b}\in T_{b}M$ be any two
\emph{tangent vectors} of $\underset{p\in U_{o}}{\bigcup}T_{p}M$.

We say that $\mathbf{v}_{a}$ is \emph{equivalent} to $\mathbf{v}_{b}$ (written
as $\mathbf{v}_{a}\mathbb{s}\mathbf{v}_{b}$) if and only if
\begin{equation}
\mathbf{v}_{a}\phi_{o}^{\mu}=\mathbf{v}_{b}\phi_{o}^{\mu},\text{ for each }%
\mu=1,\ldots,n. \label{CVS.3}%
\end{equation}

Now, since $\mathbf{v}_{a}=\mathbf{v}_{a}\phi_{o}^{\mu}\left.  \dfrac
{\partial}{\partial x_{o}^{\mu}}\right|  _{a}$ and $\mathbf{v}_{b}%
=\mathbf{v}_{b}\phi_{o}^{\mu}\left.  \dfrac{\partial}{\partial x_{o}^{\mu}%
}\right|  _{b},$ the equivalence between $\mathbf{v}_{a}$ and $\mathbf{v}_{b}$
means that the $\mu$-th \emph{contravariant components }of $\mathbf{v}_{a}$
with respect to $\{\left.  \dfrac{\partial}{\partial x_{o}^{\mu}}\right|
_{a}\}$ are equal to the $\mu$-th \emph{contravariant components }of
$\mathbf{v}_{b} $ with respect to $\{\left.  \dfrac{\partial}{\partial
x_{o}^{\mu}}\right|  _{b}\}.$

It is obvious that $\mathbb{s}$ is a well-defined \emph{ equivalence
relation}, and that it is not empty, since the tangent coordinate vectors at
any two points $a$ and $b$ belonging to $U_{o}$ are equivalent to each other,
i.e.,
\begin{equation}
\left.  \dfrac{\partial}{\partial x_{o}^{\mu}}\right\vert _{a}\mathbb{s}%
\left.  \dfrac{\partial}{\partial x_{o}^{\mu}}\right\vert _{b},\text{ for each
}\mu=1,\ldots,n. \label{CVS.4}%
\end{equation}

Let $\mathfrak{C}_{\mathbf{v}_{o}}$ be the \emph{equivalent class }of
$\mathbf{v}_{o}\in T_{o}M,$ i.e.,
\begin{equation}
\mathfrak{C}_{\mathbf{v}_{o}}=\{\mathbf{v}_{p}\text{ }/\text{ for all }p\in
U_{o}:\mathbf{v}_{p}\mathbb{s}\mathbf{v}_{o}\}. \label{CVS.5}%
\end{equation}
Let $\mathcal{U}_{o}$ be the set of \emph{all the equivalent classes} for
every $\mathbf{v}_{o}\in T_{o}M,$ i.e.,
\begin{equation}
\mathcal{U}_{o}=\{\mathfrak{C}_{\mathbf{v}_{o}}\text{ }/\text{ }\mathbf{v}%
_{o}\in T_{o}M\}. \label{CVS.6}%
\end{equation}

Such $\mathcal{U}_{o}$ has a natural structure of a \emph{real vector space}.
Indeed, such a structure is realized by defining:

The addition of vectors
\begin{equation}
\mathfrak{C}_{\mathbf{v}_{o}}\in\mathcal{U}_{o}\text{ and }\mathfrak{C}%
_{\mathbf{w}_{o}}\in\mathcal{U}_{o}\Rightarrow\mathfrak{C}_{\mathbf{v}_{o}%
}+\mathfrak{C}_{\mathbf{w}_{o}}=\mathfrak{C}_{\mathbf{v}_{o}+\mathbf{w}_{o}%
}\in\mathcal{U}_{o}. \label{CVS.7a}%
\end{equation}

The scalar multiplication of vectors by real numbers
\begin{equation}
\lambda\in\mathbb{R}\text{ and }\mathfrak{C}_{\mathbf{v}_{o}}\in
\mathcal{U}_{o}\Rightarrow\lambda\mathfrak{C}_{\mathbf{v}_{o}}=\mathfrak{C}%
_{\lambda\mathbf{v}_{o}}\in\mathcal{U}_{o}. \label{CVS.7b}%
\end{equation}

We notice that the \emph{zero vector} for $\mathcal{U}_{o}$ is given by
$0=\mathfrak{C}_{\mathbf{0}_{o}},$ where $\mathbf{0}_{o}$ is the \emph{zero
vector} for $T_{o}M$. Then, $0=\mathfrak{C}_{\mathbf{0}_{o}}$ is just the set
of all the \emph{zero tangent vectors} $\mathbf{0}_{p}\in T_{p}M$ \ for all
$p\in U_{o}$.

Now, let us take any $\mathfrak{C}_{\mathbf{v}_{o}}\in\mathcal{U}_{o}.$ Then,
by recalling that $\mathbf{v}_{o}=\mathbf{v}_{o}\phi_{o}^{\mu}\left.
\dfrac{\partial}{\partial x_{o}^{\mu}}\right|  _{o}$ we have that
\begin{equation}
\mathfrak{C}_{\mathbf{v}_{o}}=\mathbf{v}_{o}\phi_{o}^{\mu}\mathfrak{C}%
_{\left.  \frac{\partial}{\partial x_{o}^{\mu}}\right|  _{o}}. \label{CVS.8}%
\end{equation}
Eq.(\ref{CVS.8}) shows that the $n$ vectors $\mathfrak{C}_{\left.
\frac{\partial}{\partial x_{o}^{1}}\right|  _{o}},\ldots,\mathfrak{C}_{\left.
\frac{\partial}{\partial x_{o}^{n}}\right|  _{o}}$ belonging to $\mathcal{U}%
_{o}$ \emph{span} $\mathcal{U}_{o}.$

We can also check that they are in fact \emph{linearly independent}, i.e.,
\[
\lambda^{\mu}\mathfrak{C}_{\left.  \frac{\partial}{\partial x_{o}^{\mu}%
}\right\vert _{o}}=0\Rightarrow\lambda^{\mu}=0,\text{ for each }\mu
=1,\ldots,n.
\]
Indeed, by \emph{definition of equality for equivalence classes}, i.e.,
$\mathfrak{C}_{\mathbf{v}_{o}}=\mathfrak{C}_{\mathbf{w}_{o}}\Leftrightarrow
\mathbf{v}_{o}\mathbb{s}\mathbf{w}_{o},$ we have that
\[
\lambda^{\mu}\mathfrak{C}_{\left.  \frac{\partial}{\partial x_{o}^{\mu}%
}\right\vert _{o}}=0\Rightarrow\mathfrak{C}_{\lambda^{\mu}\left.
\frac{\partial}{\partial x_{o}^{\mu}}\right\vert _{o}}=\mathfrak{C}%
_{\mathbf{0}_{o}}\Rightarrow\lambda^{\mu}\left.  \frac{\partial}{\partial
x_{o}^{\mu}}\right\vert _{o}=\mathbf{0}_{o},
\]
whence, by linear independence of $\left\{  \left.  \dfrac{\partial}{\partial
x_{o}^{\mu}}\right\vert _{o}\right\}  ,$ it follows that $\lambda^{\mu}=0,$
for each $\mu=1,\ldots,n.$

We say that the above $n$ linearly independent vectors spanning $\mathcal{U}%
_{o}$ are a set of \emph{fundamental basis} vectors. We put,\emph{\ }
\begin{equation}
b_{1}=\mathfrak{C}_{\left.  \frac{\partial}{\partial x_{o}^{1}}\right|  _{o}%
},\ldots,b_{n}=\mathfrak{C}_{\left.  \frac{\partial}{\partial x_{o}^{n}%
}\right|  _{o}}. \label{CVS.9}%
\end{equation}
Then, it follows that $\dim\mathcal{U}_{o}=n$ (the dimension of $M$).

Such $\mathcal{U}_{o}$ will be called the \emph{canonical space} for the local
chart $(U_{o},\phi_{o}).$ The fundamental basis $\{b_{\mu}\}$ will be called
the\emph{\ fiducial basis} for $\mathcal{U}_{o}.$ And the real numbers
$x_{o}^{1},\ldots,x_{o}^{n}$ will be conveniently named as \emph{canonical
position coordinates} of the point $p\in U_{o}$.\footnote{If $(M,%
%TCIMACRO{\TeXButton{slg}{\slg}}%
%BeginExpansion
\slg
%EndExpansion
)$ is a $4$-dimensional Lorentzian spacetime admiting spinor fields, i.e., is
a spin manifold, then, as it is well known, it must admits a global tetrad
field (see, e.g.,\cite{rod041,moro}). In that case, the existence of the
tetrad field suggests by itself as a natural way to define an equivalence
relation between vectors at different spacetime points by the use of an
auxiliary teleparallel connection. The use of global tetrad field together
with geometrical algebra techniques \ has been used recently in an interesting
paper by Francis and Kosowsky \cite{francis}. Their results are to be compared
with the ones developed in the present series of papers.}

Our next step in order to use the geometrical and extensor calculus of
\cite{3} \ is the introduction of a fiducial Euclidean scalar product in the
canonical space $\mathcal{U}_{o}$. This is done by declaring the set
$\{b_{\mu}\}$ Euclidean orthonormal, which means, of course, that $b_{\mu
}\cdot b_{\nu}=\delta_{\mu\nu}$.

\textbf{Remark 1. }\textit{It is quite obvious that the equivalence relation
defined above is chart dependent, but this fact does not imply in any
restriction in the utilization of the methods described in this paper. Indeed,
the introduction of a canonical space }$U_{o}$\textit{ associated to a local
chart }$(U_{o},\phi_{o})$\textit{ of the given atlas of }$M$\textit{ is only a
device for the quickly application of the algebraic tools developed in
\cite{3} and has no fundamental status in the differential geometry of }%
$M$\textit{. Thus, if necessary, for the realization of some specific
calculation we simply define another canonical space associated with another
local chart }$(U_{o_{1}},\phi_{o_{1}})$\textit{ of the given atlas of }%
$M$\textit{ and use the same methodology which applies to }$(U_{o},\phi_{o}%
)$\textit{.}

\subsection{Position Vector}

The \emph{open subset} $\mathcal{U}_{o}^{\prime}\subset\mathcal{U}_{o}$,
defined by
\begin{equation}
\mathcal{U}_{o}^{\prime}=\{\lambda^{\mu}b_{\mu}\in\mathcal{U}_{o}\text{
}/\text{ }\lambda^{\mu}\in\phi_{o}^{\mu}(U_{o}),\text{ for each }\mu
=1,\ldots,n\} \label{PV.1}%
\end{equation}
will be called the \emph{position vector set} of $U_{o}.$ Of course, it is
associated to $(U_{o},\phi_{o}).$

There exists an \emph{homeomorphsim} $\iota$ between $U_{o}$ and
$\mathcal{U}_{o}^{\prime}$ which is realized by $U_{o}\ni p\mapsto\iota
_{o}(p)\in\mathcal{U}_{o}^{\prime}$ and $\mathcal{U}_{o}^{\prime}\ni
x\mapsto\iota_{o}^{-1}(x_{o})\in U_{o}$ such that
\begin{align}
\iota_{o}(p)  &  =\phi_{o}^{\mu}(p)b_{\mu},\label{PV.2a}\\
\iota_{o}^{-1}(x_{o})  &  =\phi_{o}^{-1}(b^{1}\cdot x_{o},\ldots,b^{n}\cdot
x_{o}). \label{PV.2b}%
\end{align}

As suggested by the above notations, $\iota_{o}^{-1}$ is the \emph{inverse
mapping} of $\iota_{o}.$ We have indeed that for any $p\in U_{o}$%
\begin{align*}
\iota_{o}^{-1}\circ\iota_{o}(p)  &  =\phi_{o}^{-1}(b^{1}\cdot\pi^{\mu}%
\circ\phi_{o}(p)b_{\mu},\ldots,b^{n}\cdot\pi^{\mu}\circ\phi_{o}(p)b_{\mu})\\
&  =\phi_{o}^{-1}(\pi^{1}\circ\phi_{o}(p),\ldots,\pi^{n}\circ\phi_{o}(p))\\
&  =\phi_{o}^{-1}\circ\phi_{o}(p)=p,
\end{align*}
i.e., $\iota_{o}^{-1}\circ\iota_{o}=i_{U_{o}}.$

And for any $x_{o}\in\mathcal{U}_{o}^{\prime},$
\begin{align*}
\iota_{o}\circ\iota_{o}^{-1}(x_{o})  &  =\phi_{o}^{\mu}(\phi_{o}^{-1}%
(b^{1}\cdot x_{o},\ldots,b^{n}\cdot x_{o}))b_{\mu}\\
&  =\pi^{\mu}\circ\phi_{o}\circ\phi_{o}^{-1}(b^{1}\cdot x_{o},\ldots
,b^{n}\cdot x_{o})b_{\mu}\\
&  =(b^{\mu}\cdot x_{o})b_{\mu}=x_{o},
\end{align*}
i.e., $\iota_{o}\circ\iota_{o}^{-1}=i_{\mathcal{U}_{o}^{\prime}}.$

Any \emph{point }$p\in U_{o}$ can be localized by a \emph{vector} $\iota
_{o}(p)\in\mathcal{U}_{o}^{\prime}.$ We call
\begin{equation}
x_{o}=\iota_{o}(p) \label{PV.3}%
\end{equation}
the \emph{position vector} of $p$ with respect to $(U_{o},\phi_{o}).$
Sometimes $x_{o}$ will be named as the \emph{canonical position vector} of
$p.$ By using Eq.(\ref{PV.2a}) and Eq.(\ref{CVS.1}) we can write
Eq.(\ref{PV.3}) as
\begin{equation}
x_{o}=x_{o}^{\mu}b_{\mu}. \label{PV.4}%
\end{equation}

\subsection{Canonical Algebraic Structures}

Let $\{\beta^{\mu}\}$, $\beta^{\nu}(b_{\mu})=\delta_{\mu}^{\nu}$ be the
\emph{dual basis} of $\{b_{\mu}\}.$

We \ already have equipped $\mathcal{U}_{o}$ with a fiducial Euclidean metric,
which, of course is given by $\delta_{\mu\nu}\beta^{\mu}\otimes\beta^{\nu}$.
It will be called the \emph{canonical metric}, or also, the $b$-\emph{metric,
}see \cite{11,3}

The \emph{euclidean }scalar product of $v,w\in\mathcal{U}_{o}$ corresponding
to the $b$-metric, namely $v\cdot w\in\mathbb{R}$ will be called the
\emph{canonical scalar product} \emph{of vectors}, or for short, the
$b$-\emph{scalar product.}

The $b$-\emph{reciprocal basis} of $\{b_{\mu}\},$ namely $\{b^{\mu}\}$, the
unique basis such that $b^{\nu}\cdot b_{\mu}=\delta_{\mu}^{\nu}$, coincides
with $\{b_{\mu}\},$ i.e.,
\begin{equation}
b^{\mu}=b_{\mu},\text{ for each }\mu=1,\ldots,n. \label{CAS.3b}%
\end{equation}

We denote by $\bigwedge^{k}\mathcal{U}_{o}$ ($0\leq k\leq n$) the \emph{space
of }$k$-\emph{vectors} over $\mathcal{U}_{o}$ $,$ and by $\bigwedge
\mathcal{U}_{o}$ the \emph{space of multivectors} over $\mathcal{U}_{o},$ see
\cite{2,3}. The \emph{space of} $k$-\emph{extensors} over $\mathcal{U}_{o}$
will be denoted by $k$-$ext(\bigwedge_{1}^{\diamond}\mathcal{U}_{o}%
,\ldots,\bigwedge_{k}^{\diamond}\mathcal{U}_{o};\bigwedge^{\diamond
}\mathcal{U}_{o})$ . In particular, $ext_{p}^{q}(\mathcal{U}_{o}),$
$ext(\mathcal{U}_{o})$ and $k$-$ext^{q}(\mathcal{U}_{o})$ will denote
respectively the \emph{spaces of} $(p,q)$-\emph{extensors, extensors
}and\emph{\ elementary} $k$-\emph{extensors of degree} $q$ over $\mathcal{U}%
_{o}$.

The \emph{Euclidean }scalar product of $X,Y\in\bigwedge\mathcal{U}_{o}$
relative to the \emph{euclidean metric} \emph{structure} $(\mathcal{U}%
_{o},\delta_{\mu\nu}\beta^{\mu}\otimes\beta^{\nu})$, will be denoted by
$X\cdot Y\in\mathbb{R},$ as defined in \cite{3} and will be called the
\emph{canonical scalar product of multivectors}, or for short, the
$b$-\emph{scalar product.}

The canonical algebraic structure $(\bigwedge\mathcal{U}_{o},\cdot)$ allows us
to define \emph{left} and \emph{right contracted products} of $X,Y\in
\bigwedge\mathcal{U}_{o},$ namely $X\lrcorner Y$ and $X\llcorner Y$, as
defined in \cite{3} and will be simply called the $b$-\emph{contracted
products of multivectors.} $\bigwedge\mathcal{U}_{o}$ endowed with each one of
the interior products $(\lrcorner)$ or $(\llcorner)$ is an non-associative
algebra which will be called a $b$-\emph{interior algebra of multivectors.}

The $b$-interior algebras $(\bigwedge\mathcal{U}_{o},\lrcorner)$ and
$(\bigwedge\mathcal{U}_{o},\llcorner)$ together with the exterior algebra
$(\bigwedge\mathcal{U}_{o},\wedge)$ allow us to construct a Clifford algebra
of multivectors. The \emph{Clifford product} of $X,Y\in\bigwedge
\mathcal{U}_{o},$ denoted by juxtaposition $XY\in\bigwedge\mathcal{U}_{o},$
has been defined in \cite{3}, and will be simply called the $b$-\emph{Clifford
product of multivectors}. $\bigwedge\mathcal{U}_{o}$ endowed with the
$b$-Clifford product is an associative algebra which will be called the
$b$-\emph{Clifford algebra of multivectors,} or the $b$-\emph{geometric
algebra.}

\section{Multivector Fields}

Multivector fields can be thought as sections of $%
%TCIMACRO{\dbigwedge }%
%BeginExpansion
{\displaystyle\bigwedge}
%EndExpansion
TM$, the exterior algebra bundle of multivectors. Given the equivalence
relation defined in Section 2 a multivector field $\mathbf{X\in}\sec%
%TCIMACRO{\dbigwedge }%
%BeginExpansion
{\displaystyle\bigwedge}
%EndExpansion
TU\subset\sec%
%TCIMACRO{\dbigwedge }%
%BeginExpansion
{\displaystyle\bigwedge}
%EndExpansion
TM$ ($U$ $\subset$ $U_{o}$) is, of course, represented by a mapping
\begin{equation}
X:U\rightarrow\bigwedge\mathcal{U}_{o},\label{MF.1}%
\end{equation}
which will be called a \textit{representative of the }\emph{multivector field
}$\mathbf{X}$ on $U$.\footnote{Readers, with enough knowledge of vector bundle
theory, will regonize that this is equivalent to a local trivialization of the
vector bundle $%
%TCIMACRO{\dbigwedge }%
%BeginExpansion
{\displaystyle\bigwedge}
%EndExpansion
TM$}When there is no possibility of \textit{confusion} we will call $X$ simply
a \emph{multivector field} on $U$.

The multivector function of the canonical position coordinates, namely
$X\circ\phi_{o}^{-1}$, given by
\begin{equation}
\phi_{o}(U)\ni(x_{o}^{1},\ldots,x_{o}^{n})\mapsto X\circ\phi_{o}^{-1}%
(x_{o}^{1},\ldots,x_{o}^{n})\in\bigwedge\mathcal{U}_{o}, \label{MF.2}%
\end{equation}
is called the \emph{position coordinates representation} of $X$, of course,
relative to $(U_{o},\phi_{o}),$ see \cite{8}.

The multivector function of the canonical position vector, namely $X\circ
\iota_{o}^{-1}$, given by
\begin{equation}
\iota_{o}(U)\ni x_{o}\mapsto X\circ\iota_{o}^{-1}(x_{o})\in\bigwedge
\mathcal{U}_{o}, \label{MF.3}%
\end{equation}
is called the \emph{position vector representation} of $X$ with respect to
$(U_{o},\phi_{o}),$ see \cite{9}.

In what follows we suppose that any multivector field $X$ \ used is
\emph{smooth}, i.e., $\mathcal{C}^{\infty}$ differentiable or at least enough
differentiable for our statements to hold.

The set of smooth multivector fields on $U$ will be denoted by $\mathcal{M}%
(U).$ In particular, the set of smooth scalar fields, the set of smooth vector
fields and the set of smooth $k$-vector fields ($k\geq2$) will be respectively
denoted by $\mathcal{S}(U),$ $\mathcal{V}(U)$ and $\mathcal{M}^{k}(U)$. The
\emph{identity} for $\mathcal{S}(U)$ will be denoted by $1:U\rightarrow
\mathbb{R}$ such that $1(p)=1$, and as well known, $\mathcal{M}(U)$ has a
natural structure of a\emph{\ module} over the ring (with identity)
$\mathcal{S}(U)$.

What is really important here is that $\mathcal{M}(U)$ can be endowed with
four kind of products of smooth multivector fields. Let, as in \cite{3} $\ast$
be any suitable product of multivectors either the exterior product
$(\wedge),$ the $b$-scalar product $(\cdot),$ the $b$-contracted products
$(\lrcorner,\llcorner)$ or the $b$-Clifford product. Each of these products of
multivectors induces a well-defined product of smooth multivector fields which
will be also denoted by $\ast$. The $\ast$-products of $X,Y\in\mathcal{M}(U),$
namely $X\ast Y\in\mathcal{M}(U),$ are defined by
\begin{equation}
(X\ast Y)(p)=X(p)\ast Y(p),\text{ for all }p\in U. \label{MF.6}%
\end{equation}

$\mathcal{M}(U)$ equipped with $(\wedge)$ is an associative algebra induced by
the exterior algebra of multivectors. It is called the \emph{exterior algebra
of smooth multivector fields.}

$\mathcal{M}(U)$ equipped with each of $(\lrcorner)$ or $(\llcorner)$ is a
non-associative algebra\emph{\ }induced by the respective $b$-interior algebra
of multivectors. They are called the $b$-\emph{interior algebras of smooth
multivector fields.}

$\mathcal{M}(U)$ equipped with the $b$-Clifford product is an associative
algebra induced by the $b$-Clifford algebra of multivectors. It is called the
$b$-\emph{Clifford algebra of smooth multivector fields.}

It is still possible to define in an obvious way four kind of products between
multivectors and smooth multivector fields. Indeed, let $X\in\bigwedge
\mathcal{U}_{o}$ and $Y\in\mathcal{M}(U).$ The $\ast$-product of $X$ and $Y,$
namely $X\ast Y\in\mathcal{M}(U),$ is defined by
\begin{equation}
(X\ast Y)(p)=X\ast Y(p),\text{ for all }p\in U. \label{MF.7a}%
\end{equation}

Let $X\in\mathcal{M}(U)$ and $Y\in\bigwedge\mathcal{U}_{o}$ The $*$-product of
$X$ and $Y,$ namely $X*Y\in\mathcal{M}(U),$ is defined by
\begin{equation}
(X*Y)(p)=X(p)*Y,\text{ for all }p\in U. \label{MF.7b}%
\end{equation}

\subsection{$a$-Directional Ordinary Derivative of Multivector Fields}

In \cite{8,9} we developed a complete theory of derivative operators which act
on multivector functions of real variables and on multivector functions of
multivector variables. Here, we shall need to recall only the definition of
\textit{directional derivative} of a differentiable multivector-valued
function of a vector variable, and the definition of \textit{partial
derivatives} of differentiable multivector-valued functions of several real
variables. Let $V$ be a $n$-dimensional real vector space. As usual, let us
denote by $\bigwedge V$ the \emph{space of multivectors over }$V,$ see
\cite{3}. Let $S_{1},\ldots,S_{k}$ be $k$ open subsets of $\mathbb{R}.$

Let $V\ni\mathbf{v}\mapsto F(\mathbf{v})\in\bigwedge V$ be a
\emph{differentiable }multivector function of a vector variable. Let us take
$\mathbf{a}\in V,$ the $\mathbf{a}$-directional derivative of $F$ is defined
to be
\begin{equation}
\mathbf{a\cdot}\partial_{\mathbf{v}}F(\mathbf{v})=F_{\mathbf{a}}^{\prime
}(\mathbf{v})=\underset{\lambda\rightarrow0}{\lim}\frac{F(\mathbf{v}%
+\lambda\mathbf{a})-F(\mathbf{v})}{\lambda}. \label{AOD.0a}%
\end{equation}

Let $S_{1}\times\cdots\times S_{k}\ni(\lambda^{1},\ldots,\lambda^{k})\mapsto
f(\lambda^{1},\ldots,\lambda^{k})\in\bigwedge V$ be a \emph{differentiable}
multivector function of $k$ real variables. The $\lambda^{j}$-partial
derivative of $f$ (with $1\leq j\leq k$) is defined to be
\begin{align}
\frac{\partial f}{\partial\lambda^{j}}(\lambda^{1},\ldots,\lambda^{k})  &
=f^{\prime(j)}(\lambda^{1},\ldots,\lambda^{k})\nonumber\\
&  =\underset{\mu\rightarrow0}{\lim}\frac{f(\lambda^{1},\ldots,\lambda^{j}%
+\mu,\ldots,\lambda^{k})-f(\lambda^{1},\ldots,\lambda^{k})}{\mu}
\label{AOD.0b}%
\end{align}

\noindent\textbf{Proposition 1.}\textit{ Let us take} $a\in\mathcal{U}_{o}$.
\textit{For any smooth multivector field} $X,$ \textit{the} $a$%
-\textit{directional ordinary derivative of the position vector representation
of} $X$ \textit{with respect to }$(U_{o},\phi_{o}),$\textit{ namely }%
$a\cdot\partial_{x_{o}}X\circ\iota_{o}^{-1},$\textit{ is related to the
}$x_{o}^{\mu}$\textit{-partial derivatives of the position}
\textit{coordinates representation of }$X$\textit{ with respect to }%
$(U_{o},\phi_{o}),$\textit{ namely }$\dfrac{\partial}{\partial x_{o}^{\mu}%
}X\circ\phi_{o}^{-1},$\textit{ by the identity }%
\begin{equation}
(a\cdot\partial_{x_{o}}X\circ\iota_{o}^{-1})\circ\iota_{o}=a\cdot b^{\mu
}(\dfrac{\partial}{\partial x_{o}^{\mu}}X\circ\phi_{o}^{-1})\circ\phi_{o}.
\label{AOD.1}%
\end{equation}

\noindent\textbf{Proof}

It is enough to verify that
\[
(b_{\mu}\cdot\partial_{x_{o}}X\circ\iota_{o}^{-1})\circ\iota_{o}%
=(\dfrac{\partial}{\partial x_{o}^{\mu}}X\circ\phi_{o}^{-1})\circ\phi
_{o},\text{ for each }\mu=1,\ldots,n.
\]

Since the position vector $x_{o}$ is an interior point of $\iota_{o}(U),$
there is some $\varepsilon$neighborhood, say $\mathcal{N}_{x_{o}}%
(\varepsilon),$ such that $\mathcal{N}_{x_{o}}(\varepsilon)\subseteq\iota
_{o}(U).$ Now, if we take $\lambda\in\mathbb{R}$ such that $0<\left\vert
\lambda\right\vert <\varepsilon,$ it follows that $x_{o}+\lambda b_{\mu}%
\in\mathcal{N}_{x_{o}}(\varepsilon).$ It follows that $x_{o}+\lambda b_{\mu
}\in\iota_{o}(U)$, and there exist $\iota_{o}^{-1}(x_{o}+\lambda b_{\mu})\in
U$ and $X\circ\iota_{o}^{-1}(x_{o}+\lambda b_{\mu})\in\bigwedge\mathcal{U}%
_{o}.$

Hence, by taking into account Eq.(\ref{PV.2b}), we get the following identity
\begin{align*}
&  \frac{X\circ\iota_{o}^{-1}(x_{o}+\lambda b_{\mu})-X\circ\iota_{o}%
^{-1}(x_{o})}{\lambda}\\
&  =\frac{X\circ\phi_{o}^{-1}(b^{1}\cdot x_{o}+\lambda\delta_{\mu}^{1}%
,\ldots,b^{n}\cdot x_{o}+\lambda\delta_{\mu}^{n})-X\circ\phi_{o}^{-1}%
(b^{1}\cdot x_{o},\ldots,b^{n}\cdot x_{o})}{\lambda}.
\end{align*}

Now, by taking limits for $\lambda\rightarrow0$ on these multivector functions
of the real variable $\lambda,$ and using Eq.(\ref{PV.2b}) once again, we have
indeed that
\begin{align*}
b_{\mu}\cdot\partial_{x_{o}}X\circ\iota_{o}^{-1}(x_{o})  &  =\delta_{\mu}%
^{\nu}\frac{\partial}{\partial x_{o}^{\nu}}X\circ\phi_{o}^{-1}(b^{1}\cdot
x_{o},\ldots,b^{n}\cdot x_{o})\\
&  =(\frac{\partial}{\partial x_{o}^{\mu}}X\circ\phi_{o}^{-1})\circ\phi
_{o}\circ\iota_{o}^{-1}(x_{o}),
\end{align*}
and the required result follows.$\blacksquare$

Let $X$ be a smooth multivector field on $U,$ and let us take $a\in
\mathcal{U}_{o}.$ The smooth multivector field on $U,$ namely $a\cdot
\partial_{o}X,$ defined as $U\ni p\mapsto$ $a\cdot\partial_{o}X(p)\in
\bigwedge\mathcal{U}_{o}$ such that
\begin{equation}
a\cdot\partial_{o}X(p)=(a\cdot\partial_{x_{o}}X\circ\iota_{o}^{-1})(x_{o}),
\label{AOD.2}%
\end{equation}
where $x_{o}=\iota_{o}(p),$ will be called the \emph{canonical} $a$%
-\emph{directional ordinary derivative} ($a$-\emph{DOD}) of $X.$ Sometimes,
$a\cdot\partial_{o}$ will be named as the \emph{canonical a-directional
ordinary derivative operator} ($a$-\emph{DODO}).

We see that the position vector representation of $a\cdot\partial_{o}X$ is
equal to the $a$-directional derivative of the position vector representation
of $X,$ of course, both of them respect to $(U_{o},\phi_{o}). $

Taking into account Eq.(\ref{AOD.1}), we have for any $a\in\mathcal{U}_{o},$
and for all $X\in\mathcal{M}(U)$ that
\begin{equation}
a\cdot\partial_{o}X(p)=a\cdot\partial_{x_{o}}X\circ\iota_{o}^{-1}%
(x_{o})=a\cdot b^{\mu}\dfrac{\partial}{\partial x_{o}^{\mu}}X\circ\phi
_{o}^{-1}(x_{o}^{1},\ldots,x_{o}^{n}), \label{AOD.3}%
\end{equation}
where $x_{o}=\iota_{o}(p)$ and $(x_{o}^{1},\ldots,x_{o}^{n})=\phi_{o}(p)$.

It is convenient and useful to generalize the notion of canonical
$a$-\emph{DOD} of a smooth multivector field whenever $a$ is any smooth vector field.

Let us take $a\in\mathcal{V}(U).$ The canonical $a$-\emph{DOD} of $X,$ also
denoted by $a\cdot\partial_{o}X,$ is defined to be $U\ni p\mapsto
a\cdot\partial_{o}X(p)\in\bigwedge\mathcal{U}_{o}$ such that
\begin{equation}
a\cdot\partial_{o}X(p)=a\circ\iota_{o}^{-1}(x_{o})\cdot\partial_{x_{o}}%
X\circ\iota_{o}^{-1}(x_{o}), \label{AOD.3a}%
\end{equation}
where $x_{o}=\iota_{o}(p).$ In agreement with Eq.(\ref{AOD.1}), we may also
write
\begin{equation}
a\cdot\partial_{o}X(p)=a\circ\phi_{o}^{-1}(x_{o}^{1},\ldots,x_{o}^{n})\cdot
b^{\mu}\dfrac{\partial}{\partial x_{o}^{\mu}}X\circ\phi_{o}^{-1}(x_{o}%
^{1},\ldots,x_{o}^{n}), \label{AOD.3b}%
\end{equation}
where $(x_{o}^{1},\ldots,x_{o}^{n})=\phi_{o}(p).$

We summarize now the basic properties which are satisfied by the
$a$-\emph{DOD's.\vspace{0.1in}}

\textbf{i.} For any $a\in\mathcal{U}_{o}$ or $a\in\mathcal{V}(U)$ the
canonical $a$-\emph{DODO}, namely $a\cdot\partial_{o},$ is grade-preserving,
i.e.,
\begin{equation}
\text{if }X\in\mathcal{M}^{k}(U),\text{ then }a\cdot\partial_{o}%
X\in\mathcal{M}^{k}(U). \label{AOD.4}%
\end{equation}

\textbf{ii.} For any $a,a^{\prime}\in\mathcal{U}_{o}$ or $a,a^{\prime}%
\in\mathcal{V}(U),$ and $\alpha,\alpha^{\prime}\in\mathbb{R},$ and for all
$X\in\mathcal{M}(U),$ it holds
\begin{equation}
(\alpha a+\alpha^{\prime}a^{\prime})\cdot\partial_{o}X=\alpha a\cdot
\partial_{o}X+\alpha^{\prime}a^{\prime}\cdot\partial_{o}X, \label{AOD.5}%
\end{equation}

\textbf{iii. }For any $a,a^{\prime}\in\mathcal{U}_{o}$ or $a,a^{\prime}%
\in\mathcal{V}(U),$ and $f,f^{\prime}\in\mathcal{S}(U),$ and for all
$X\in\mathcal{M}(U),$ it holds
\begin{equation}
(fa+f^{\prime}a^{\prime})\cdot\partial_{o}X=fa\cdot\partial_{o}X+f^{\prime
}a^{\prime}\cdot\partial_{o}X, \label{AOD.5a}%
\end{equation}

\textbf{iv.} For any $a\in\mathcal{U}_{o}$ or $a\in\mathcal{V}(U),$ and for
all $f\in\mathcal{S}(U)$ and $X,Y\in\mathcal{M}(U),$ it holds
\begin{align}
a\cdot\partial_{o}(X+Y)  &  =a\cdot\partial_{o}X+a\cdot\partial_{o}%
Y\label{AOD.6a}\\
a\cdot\partial_{o}(fX)  &  =(a\cdot\partial_{o}f)X+f(a\cdot\partial_{o}X).
\label{AOD.6b}%
\end{align}

\textbf{v.} Let $\ast$ mean any multivector product either $(\wedge),$
$(\cdot)$ $(\lrcorner,\llcorner)$ or $(b$-\emph{Clifford product}$).$ For any
$a\in\mathcal{U}_{o}$ or $a\in\mathcal{V}(U),$ and for all $X,Y\in
\mathcal{M}(U),$ it holds
\begin{equation}
a\cdot\partial_{o}(X\ast Y)=(a\cdot\partial_{o}X)\ast Y+X\ast(a\cdot
\partial_{o}Y). \label{AOD.7}%
\end{equation}
It is a Leibniz-like rule for any suitable multivector product as introduced
in Section 3.

\textbf{Remark 2. }\textit{The concept of }$a$\textit{- direction ordinary
derivative is obviously chart dependent. It has been introduced as part of a
calculation device, whose main utility results from the fact that }%
$a\cdot\partial_{o}$\textit{ is a connection on }$U_{o}$\textit{, a crucial
fact that will be used in the following papers of the series.}

\subsection{Lie Algebra of Smooth Vector Fields}

The module $\mathcal{V}(U)$ of the smooth vector fields on $U$, can be
equipped with a well-defined \emph{Lie product}. Here we briefly recall how to
write this concept and their properties in the present approach. Let
$a,b\in\mathcal{V}(U)$, the (canonical) Lie product $a$ and $b,$ namely
$[a,b],$ is defined by
\begin{equation}
\lbrack a,b]=a\cdot\partial_{o}b-b\cdot\partial_{o}a. \label{LA.1}%
\end{equation}
By using Eq.(\ref{AOD.3a}) and Eq.(\ref{AOD.3b}), we can get two noticeable
formulas for this Lie product. One in terms of the (vector)-directional
derivative operators with respect to $x_{o}$ acting on the canonical position
vector representation of $a$ and $b,$ i.e.,
\begin{equation}
\lbrack a,b](p)=a\circ\iota_{o}^{-1}(x_{o})\cdot\partial_{x_{o}}b\circ
\iota_{o}^{-1}(x_{o})-b\circ\iota_{o}^{-1}(x_{o})\cdot\partial_{x_{o}}%
a\circ\iota_{o}^{-1}(x_{o}) \label{LA.1a}%
\end{equation}
And, another formula involving the fiducial basis $\{b_{\mu}\}$, the
$x_{o}^{\mu}$-partial derivative operators $\dfrac{\partial}{\partial
x_{o}^{\mu}}$ and the canonical position coordinates representation of $a$ and
$b,$ i.e.,
\begin{align}
\lbrack a,b](p)  &  =a\circ\phi_{o}^{-1}(x_{o}^{1},\ldots,x_{o}^{n})\cdot
b^{\mu}\frac{\partial}{\partial x_{o}^{\mu}}b\circ\phi_{o}^{-1}(x_{o}%
^{1},\ldots,x_{o}^{n})\nonumber\\
&  -b\circ\phi_{o}^{-1}(x_{o}^{1},\ldots,x_{o}^{n})\cdot b^{\mu}\frac
{\partial}{\partial x_{o}^{\mu}}a\circ\phi_{o}^{-1}(x_{o}^{1},\ldots,x_{o}%
^{n}). \label{LA.1b}%
\end{align}

We summarize the basic properties satisfied by the Lie product.

\textbf{i.} For all $a,a^{\prime},b,b^{\prime}\in\mathcal{V}(U)$
\begin{align}
\lbrack a+a^{\prime},b]  &  =[a,b]+[a^{\prime},b],\label{LA.2a}\\
\lbrack a,b+b^{\prime}]  &  =[a,b]+[a,b^{\prime}]\text{ (distributive laws).}
\label{LA.2b}%
\end{align}

\textbf{ii. }For all $f\in\mathcal{S}(U),$ and $a,b\in\mathcal{V}(U)$
\begin{align}
\lbrack fa,b]  &  =f[a,b]-(b\cdot\partial_{o}f)a,\label{LA.3a}\\
\lbrack a,fb]  &  =(a\cdot\partial_{o}f)b+f[a,b]. \label{LA.3b}%
\end{align}

\textbf{iii. F}or any $a,b\in\mathcal{V}(U),$ and for all $X\in\mathcal{M}%
(U),$ it holds\footnote{The \emph{commutator} of $a\cdot\partial_{o}$ and
$b\cdot\partial_{o},$ namely $[a\cdot\partial_{o},b\cdot\partial_{o}],$ is
defined by $[a\cdot\partial_{o},b\cdot\partial_{o}]X=a\cdot\partial_{o}%
(b\cdot\partial_{o}X)-b\cdot\partial_{o}(a\cdot\partial_{o}X),$ for all
$X\in\mathcal{M}(U)$.}
\begin{equation}
\lbrack a\cdot\partial_{o},b\cdot\partial_{o}]X=[a,b]\cdot\partial_{o}X.
\label{LA.4}%
\end{equation}

The proof of this result is as follows. By using Eqs.(\ref{AOD.6a}) and
(\ref{AOD.6b}), and Eq.(\ref{AOD.7}), we have
\begin{align}
a\cdot\partial_{o}(b\cdot\partial_{o}X)  &  =a\cdot b^{\mu}b_{\mu}%
\cdot\partial_{o}(b\cdot b^{\nu}b_{\nu}\cdot\partial_{o}X)\nonumber\\
&  =a\cdot b^{\mu}(b_{\mu}\cdot\partial_{o}b)\cdot b^{\nu}b_{\nu}\cdot
\partial_{o}X+a\cdot b^{\mu}b\cdot b^{\nu}b_{\mu}\cdot\partial_{o}(b_{\nu
}\cdot\partial_{o}X),\nonumber\\
&  =(a\cdot\partial_{o}b)\cdot\partial_{o}X+a\cdot b^{\mu}b\cdot b^{\nu}%
b_{\mu}\cdot\partial_{o}(b_{\nu}\cdot\partial_{o}X). \label{LA.4a}%
\end{align}
And, by interchanging the letters $a$ and $b,$ and re-naming indices, we have
\begin{align}
b\cdot\partial_{o}(a\cdot\partial_{o}X)  &  =(b\cdot\partial_{o}%
a)\cdot\partial_{o}X+b\cdot b^{\mu}a\cdot b^{\nu}b_{\mu}\cdot\partial
_{o}(b_{\nu}\cdot\partial_{o}X),\nonumber\\
&  =(b\cdot\partial_{o}a)\cdot\partial_{o}X+a\cdot b^{\mu}b\cdot b^{\nu}%
b_{\nu}\cdot\partial_{o}(b_{\mu}\cdot\partial_{o}X). \label{LA.4b}%
\end{align}

Now, subtracting Eq.(\ref{LA.4b}) from Eq.(\ref{LA.4a}), we get
\[
\lbrack a\cdot\partial_{o},b\cdot\partial_{o}]X=[a,b]\cdot\partial_{o}X+a\cdot
b^{\mu}b\cdot b^{\nu}(b_{\mu}\cdot\partial_{o}(b_{\nu}\cdot\partial
_{o}X)-b_{\nu}\cdot\partial_{o}(b_{\mu}\cdot\partial_{o}X)).
\]
Then, by recalling the obvious property $b_{\mu}\cdot\partial_{o}(b_{\nu}%
\cdot\partial_{o}X)=b_{\nu}\cdot\partial_{o}(b_{\mu}\cdot\partial_{o}X),$ we
finally get the required result.

\textbf{iv.} For all $a,b,c\in\mathcal{V}(U)$
\begin{equation}
\lbrack a,[b,c]]+[b,[c,a]]+[c,[a,b]]=0, \label{LA.5}%
\end{equation}
called \emph{Jacobi's identity.}

$\mathcal{M}(U)$ endowed with the Lie product is (for each $p\in U$) a
\emph{Lie algebra} which will be called the \emph{Lie algebra of smooth vector
fields.}

\subsection{Hestenes Derivatives}

Let us take any two \emph{reciprocal frame fields}\footnote{A \emph{frame
field} on $U_{o}$ is a set of $n$ smooth vector fields $e_{1},\ldots,e_{n}%
\in\mathcal{V}(U)$ such that for each $p\in U$ the set of the $n$ vectors
$e_{1}(p),\ldots,e_{n}(p)\in\mathcal{U}_{o}$ is a \emph{basis} for
$\mathcal{U}_{o}.$ In particular, $\{b_{\mu}\}$ can be taken as a
\emph{constant frame field} on $U_{o},$ i.e., $a\cdot\partial_{o}b_{\mu}=0,$
for each $\mu=1,\ldots,n.$ It is called the \emph{fiducial frame field}.} on
$U,$ namely $\{e_{\mu}\}$ and $\{e^{\mu}\},$ i.e., $e_{\mu}\cdot e^{\nu
}=\delta_{\mu}^{\nu}.$ Let $X$ be a smooth multivector field on $U$. We can
introduce exactly three smooth derivative-like multivector fields on $U,$
namely $\partial_{o}\ast X,$ defined by
\begin{equation}
\partial_{o}\ast X=e^{\mu}\ast e_{\mu}\cdot\partial_{o}X=e_{\mu}\ast e^{\mu
}\cdot\partial_{o}X, \label{HD.1}%
\end{equation}
where $\ast$ means any multivector product $(\wedge)$, $(\lrcorner)$ or
$(b$-\emph{Clifford product}$).$

Note that a smooth multivector field defined by Eq.(\ref{HD.1}) does not
depend on the choice of $\{e_{\mu}\}$ and $\{e^{\mu}\}.$

By recalling Eq.(\ref{AOD.5a}) and using Eq.(\ref{AOD.3a}) we can get a
remarkable formula for calculating $\partial_{o}\ast X$ which involves only
the fiducial frame field $\{b_{\mu}\}$ and the $b_{\mu}$-directional operators
$b_{\mu}\cdot\partial_{x_{o}}$ acting on the canonical position vector
representation of $X,$ i.e.,
\begin{equation}
\partial_{o}\ast X(p)=b^{\mu}\ast b_{\mu}\cdot\partial_{x_{o}}X\circ\iota
_{o}^{-1}(x_{o}), \label{HD.1a}%
\end{equation}
where $x_{o}=\iota_{o}(p).$

By using Eq.(\ref{AOD.3b}) we can get another formula for calculating
$\partial_{o}\ast X$ which involves only the fiducial frame fields $\{b_{\mu
}\}$ and the $x_{o}^{\mu}$-partial derivative operators $\dfrac{\partial
}{\partial x_{o}^{\mu}}$ acting on the canonical position coordinate
representation of $X,$ i.e.,
\begin{equation}
\partial_{o}\ast X(p)=b^{\mu}\ast\dfrac{\partial}{\partial x_{o}^{\mu}}%
X\circ\phi_{o}^{-1}(x_{o}^{1},\ldots,x_{o}^{n}), \label{HD.1b}%
\end{equation}
where $(x_{o}^{1},\ldots,x_{o}^{n})=\phi_{o}(p).$

We will call $\partial_{o}\wedge X,$ $\partial_{o}\lrcorner X$ and
$\partial_{o}X$ (i.e., $\ast\equiv b$-\emph{Clifford product}) respectively
the \emph{curl}, \emph{left} $b$-\emph{contracted divergence} and
$b$-\emph{gradient} of $X.$ An interesting and useful relationship relating
these derivatives is
\begin{equation}
\partial_{o}X=\partial_{o}\lrcorner X+\partial_{o}\wedge X. \label{HD.2}%
\end{equation}
In what follows we call these derivatives by \emph{Hestenes derivatives.}

We end this section presenting three remarkable multivector identities
involving the Hestenes derivatives.%

\begin{align}
(\partial_{o}\wedge X)\cdot Y+X\cdot(\partial_{o}\lrcorner Y)  &
=\partial_{o}\cdot(\partial_{n}(n\wedge X)\cdot Y),\label{DH.3a}\\
(\partial_{o}\lrcorner X)\cdot Y+X\cdot(\partial_{o}\wedge Y)  &
=\partial_{o}\cdot(\partial_{n}(n\lrcorner X)\cdot Y),\label{DH.3b}\\
(\partial_{o}X)\cdot Y+X\cdot(\partial_{o}Y)  &  =\partial_{o}\cdot
(\partial_{n}(nX)\cdot Y). \label{DH.3c}%
\end{align}
They are used in the \emph{Lagrangian theory of multivector fields.}%
\footnote{A preliminary construction of that theory on Minkowski spacetime has
been given in \cite{10}, and with more details in \cite{rodoliv2006}.}

\section{Extensor Fields}

Extensor fields are sections of an appropriate vector bundle $E(%
%TCIMACRO{\dbigwedge }%
%BeginExpansion
{\displaystyle\bigwedge}
%EndExpansion
TM)$ (which will not be described here). What is important for us is that
given the equivalence relation defined in Section 2 any $\mathbf{t\in}\sec E(%
%TCIMACRO{\dbigwedge }%
%BeginExpansion
{\displaystyle\bigwedge}
%EndExpansion
TU)\subset E(%
%TCIMACRO{\dbigwedge }%
%BeginExpansion
{\displaystyle\bigwedge}
%EndExpansion
TM)$, with $U$ $\subset$ $U_{o}$ has a representative given by the
mapping\footnote{As in the case of representatives of multivector fields, this
constructuion is equivalent to a local trivialization of $E(%
%TCIMACRO{\dbigwedge }%
%BeginExpansion
{\displaystyle\bigwedge}
%EndExpansion
TM)$.}
\begin{equation}
t:U\rightarrow k\text{-}ext(%
%TCIMACRO{\dbigwedge \nolimits_{1}^{\diamond}}%
%BeginExpansion
{\displaystyle\bigwedge\nolimits_{1}^{\diamond}}
%EndExpansion
\mathcal{U}_{o},\ldots,%
%TCIMACRO{\dbigwedge \nolimits_{k}^{\diamond}}%
%BeginExpansion
{\displaystyle\bigwedge\nolimits_{k}^{\diamond}}
%EndExpansion
\mathcal{U}_{o};%
%TCIMACRO{\dbigwedge \nolimits^{\diamond}}%
%BeginExpansion
{\displaystyle\bigwedge\nolimits^{\diamond}}
%EndExpansion
\mathcal{U}_{o}), \label{EF.1}%
\end{equation}
called a representative of $k$-\emph{extensor field} on $U$, or when no
confusion arises, simply by $k$-\emph{extensor field} on $U$. Then, for each
$p\in U:t_{(p)}$ is a $k$-linear mapping from $\bigwedge_{1}^{\diamond
}\mathcal{U}_{o}\times\ldots\times\bigwedge_{k}^{\diamond}\mathcal{U}_{o}$ to
$\bigwedge^{\diamond}\mathcal{U}_{o}$, $t_{(p)}:%
%TCIMACRO{\dbigwedge \nolimits_{1}^{\diamond}}%
%BeginExpansion
{\displaystyle\bigwedge\nolimits_{1}^{\diamond}}
%EndExpansion
\mathcal{U}_{o},\ldots,%
%TCIMACRO{\dbigwedge \nolimits_{k}^{\diamond}}%
%BeginExpansion
{\displaystyle\bigwedge\nolimits_{k}^{\diamond}}
%EndExpansion
\mathcal{U}_{o};%
%TCIMACRO{\dbigwedge \nolimits^{\diamond}}%
%BeginExpansion
{\displaystyle\bigwedge\nolimits^{\diamond}}
%EndExpansion
\mathcal{U}_{o}$.

The $k$-extensor function of the position coordinates $(x_{o}^{1},\ldots
,x_{o}^{n}),$ namely $t_{\circ\phi_{o}^{-1}}$, given by
\begin{equation}
\phi_{o}(U)\ni(x_{o}^{1},\ldots,x_{o}^{n})\mapsto t_{\circ\phi_{o}^{-1}%
(x_{o}^{1},\ldots,x_{o}^{n})}\in k\text{-}ext(%
%TCIMACRO{\dbigwedge \nolimits_{1}^{\diamond}}%
%BeginExpansion
{\displaystyle\bigwedge\nolimits_{1}^{\diamond}}
%EndExpansion
\mathcal{U}_{o},\ldots,%
%TCIMACRO{\dbigwedge \nolimits_{k}^{\diamond}}%
%BeginExpansion
{\displaystyle\bigwedge\nolimits_{k}^{\diamond}}
%EndExpansion
\mathcal{U}_{o};%
%TCIMACRO{\dbigwedge \nolimits^{\diamond}}%
%BeginExpansion
{\displaystyle\bigwedge\nolimits^{\diamond}}
%EndExpansion
\mathcal{U}_{o}), \label{EF.2}%
\end{equation}
is called the \emph{position coordinates representation} of $t$ with respect
to $(U_{o},\phi_{o})$.

The $k$-extensor function of the position vector $x_{o},$ namely
$t_{\circ\iota_{o}^{-1}},$ given by
\begin{equation}
\iota_{o}(U)\ni x_{o}\mapsto t_{\circ\iota_{o}^{-1}(x_{o})}\in k\text{-}ext(%
%TCIMACRO{\dbigwedge \nolimits_{1}^{\diamond}}%
%BeginExpansion
{\displaystyle\bigwedge\nolimits_{1}^{\diamond}}
%EndExpansion
\mathcal{U}_{o},\ldots,%
%TCIMACRO{\dbigwedge \nolimits_{k}^{\diamond}}%
%BeginExpansion
{\displaystyle\bigwedge\nolimits_{k}^{\diamond}}
%EndExpansion
\mathcal{U}_{o};%
%TCIMACRO{\dbigwedge \nolimits^{\diamond}}%
%BeginExpansion
{\displaystyle\bigwedge\nolimits^{\diamond}}
%EndExpansion
\mathcal{U}_{o}), \label{EF.3}%
\end{equation}
is called the\emph{\ position vector representation }of $t$ with respect to
$(U_{o},\phi_{o}).$

All $k$-extensor field $t$ used in what follows are supposed smooth and we can
easily prove the following proposition.

\textbf{Proposition 2.} \textit{A }$k$\textit{-extensor field }$t$\textit{ is
smooth if and only if the multivector field defined as }$U\ni p\mapsto
t_{(p)}(X_{1}(p),\ldots,X_{k}(p))\in\bigwedge^{\diamond}U_{o}$\textit{ is
itself smooth, for all }$X_{1}\in M_{1}^{\diamond}(U),\ldots,X_{k}\in
M_{k}^{\diamond}(U)$\textit{.}

Associated to any smooth $k$-extensor field $t$ we can define a \emph{linear
operator}, namely $t_{op},$ which takes $k$-uples of smooth multivector fields
belonging to $\mathcal{M}_{1}^{\diamond}(U)\times\ldots\times\mathcal{M}%
_{k}^{\diamond}(U)$ into smooth multivector fields belonging to $\mathcal{M}%
^{\diamond}(U),$ i.e.,
\begin{equation}
t_{op}(X_{1},\ldots,X_{k})(p)=t_{(p)}(X_{1}(p),\ldots,X_{k}(p)),\text{ for
each }p\in U. \label{EF.4}%
\end{equation}

Reciprocally, given a linear operator $t_{op}$ we can prove (without
difficulty) that there is an \emph{unique} smooth $k$-extensor field $t$ such
that Eq.(\ref{EF.4}) holds.

The set of smooth $k$-extensor fields on $U$ has a natural structure of
\emph{module} over the ring (with identity) $\mathcal{S}(U).$ It will be
denoted in what follows by the suggestive notation $k$-$ext(\mathcal{M}%
_{1}^{\diamond}(U),\ldots,\mathcal{M}_{k}^{\diamond}(U);\mathcal{M}^{\diamond
}(U)).$

\subsection{$a$-Directional Ordinary Derivative of Extensor Fields}

Let $t$ be a smooth $k$-extensor field on $U,$ and let us take $a\in
\mathcal{U}_{o}.$ The smooth $k$-extensor field on $U,$ namely $a\cdot
\partial_{o}t,$ such that for all smooth multivector fields $X_{1}%
\in\mathcal{M}_{1}^{\diamond}(U),\ldots,X_{k}\in\mathcal{M}_{k}^{\diamond
}(U)$
\begin{align}
(a\cdot\partial_{o}t)_{(p)}(X_{1}(p),\ldots,X_{k}(p))  &  =a\cdot\partial
_{o}(t_{(p)}(\ldots))-t_{(p)}(a\cdot\partial_{o}X_{1}(p),\ldots)\nonumber\\
&  -\cdots-t_{(p)}(\ldots,a\cdot\partial_{o}X_{k}(p)) \label{ODE.1}%
\end{align}
for each $p\in U,$ will be called the \emph{canonical }$a$-\emph{directional
ordinary derivative} ($a$-\emph{DOD}) of $t.$

Note that on the right side of Eq.(\ref{ODE.1}) there just appears the
$a$-\emph{DOD }of the smooth multivector field $U\ni p\mapsto t_{(p)}%
(X_{1}(p),\ldots,X_{k}(p))\in\bigwedge^{\diamond}\mathcal{U}_{o}.$ As usual,
when no confusion arises, we will write the definition given by
Eq.(\ref{ODE.1}) by omitting the letter $p.$

We note that the algebraic object $a\cdot\partial_{o}t$ as defined by
Eq.(\ref{ODE.1}) is indeed a smooth $k$-extensor field. The $k$-extensor
character and the smoothness for it follow directly from the respective
properties of $t.$ Also, note that in definition given by Eq.(\ref{ODE.1}),
$a$ could be a smooth vector field on $U$ instead a vector of $\mathcal{U}_{o}
$.

We present now some of the most basic properties satisfied by the
$a$-\emph{DOD's }of smooth $k$-extensor fields.\vspace{0.1in}

\textbf{i.} For any of either $a\in\mathcal{U}_{o}$ or $a\in\mathcal{V}(U)$
the canonical $a$-\emph{DODO,} namely $a\cdot\partial_{o}$ (whenever it is
acting on smooth $k$-extensor fields), preserves the extensor type, i.e.,
\begin{align}
\text{if }t  &  \in k\text{-}ext(\mathcal{M}_{1}^{\diamond}(U),\ldots
,\mathcal{M}_{k}^{\diamond}(U);\mathcal{M}^{\diamond}(U)),\nonumber\\
\text{then }a\cdot\partial_{o}t  &  \in k\text{-}ext(\mathcal{M}_{1}%
^{\diamond}(U),\ldots,\mathcal{M}_{k}^{\diamond}(U);\mathcal{M}^{\diamond
}(U)). \label{ODE.2}%
\end{align}

\textbf{ii.} For any of either $a,a^{\prime}\in\mathcal{U}_{o}$ or
$a,a^{\prime}\in\mathcal{V}(U),$ and $\alpha,\alpha^{\prime}\in\mathbb{R},$
and for all $t\in k$-$ext(\mathcal{M}_{1}^{\diamond}(U),\ldots,\mathcal{M}%
_{k}^{\diamond}(U);\mathcal{M}^{\diamond}(U)),$ it holds
\begin{equation}
(\alpha a+\alpha^{\prime}a^{\prime})\cdot\partial_{o}t=\alpha a\cdot
\partial_{o}t+\alpha^{\prime}a^{\prime}\cdot\partial_{o}t. \label{ODE.3}%
\end{equation}

\textbf{iii. }For any of either $a,a^{\prime}\in\mathcal{U}_{o}$ or
$a,a^{\prime}\in\mathcal{V}(U),$ and $f,f^{\prime}\in\emph{S}(U),$ and for all
$t\in k$-$ext(\mathcal{M}_{1}^{\diamond}(U),\ldots,\mathcal{M}_{k}^{\diamond
}(U);\mathcal{M}^{\diamond}(U)),$ it holds
\begin{equation}
(fa+f^{\prime}a^{\prime})\cdot\partial_{o}t=fa\cdot\partial_{o}t+f^{\prime
}a^{\prime}\cdot\partial_{o}t. \label{ODE.3a}%
\end{equation}

\textbf{iv.} For any of either $a\in\mathcal{U}_{o}$ or $a\in\mathcal{V}(U),$
and for all $f\in\mathcal{S}(U),$ and for all $t,u\in k$-$ext(\mathcal{M}%
_{1}^{\diamond}(U),\ldots,\mathcal{M}_{k}^{\diamond}(U);\mathcal{M}^{\diamond
}(U)),$ it holds
\begin{align}
a\cdot\partial_{o}(t+u)  &  =a\cdot\partial_{o}t+a\cdot\partial_{o}%
u,\label{ODE.4a}\\
a\cdot\partial_{o}(ft)  &  =(a\cdot\partial_{o}f)t+f(a\cdot\partial_{o}t).
\label{ODE.4b}%
\end{align}

\textbf{iv.} For any of either $a\in\mathcal{U}_{o}$ or $a\in\mathcal{V}(U),$
and for all $t\in k$-$ext(\mathcal{M}_{1}^{\diamond}(U);\mathcal{M}^{\diamond
}(U)),$ it holds
\begin{equation}
(a\cdot\partial_{o}t)^{\dagger}=a\cdot\partial_{o}t^{\dagger}. \label{ODE.5}%
\end{equation}

To prove Eq.(\ref{ODE.5}) let us take $X_{1}\in\mathcal{M}_{1}^{\diamond}(U)$
and $X\in\mathcal{M}^{\diamond}(U).$ By recalling the fundamental property of
the adjoint operator $\left.  {}\right.  ^{\dagger},$ and using
Eq.(\ref{AOD.7}), we can write
\begin{align*}
(a\cdot\partial_{o}t)^{\dagger}(X)\cdot X_{1}  &  =(a\cdot\partial_{o}%
t)(X_{1})\cdot X\\
&  =a\cdot\partial_{o}(t(X_{1})\cdot X)-t(X_{1})\cdot a\cdot\partial
_{o}X-t(a\cdot\partial_{o}X_{1})\cdot X\\
&  =a\cdot\partial_{o}(t^{\dagger}(X)\cdot X_{1})-t^{\dagger}(a\cdot
\partial_{o}X)\cdot X_{1}-t^{\dagger}(X)\cdot a\cdot\partial_{o}X_{1},\\
&  =(a\cdot\partial_{o}t^{\dagger})(X)\cdot X_{1}.
\end{align*}
Hence, by the non-degeneracy of scalar product, the expected result
immediately follows.

\section{How to Use Coordinate Systems Different from $\{x^{\mu}\}$}

Let $U$ be an \emph{open subset} of $U_{o}$, and let $(U,\phi)$ be a
\emph{local coordinate system} on $U$ \emph{compatible} with $(U_{o},\phi
_{o})$. This implies that any \emph{point} $p\in U$ can indeed be localized by
a $n$-\emph{uple of} \emph{real numbers} $\phi(p)\in\mathbb{R}^{n}.$

The $n$ \emph{scalar fields} on $U$
\begin{equation}
\phi^{\mu}:U\rightarrow\mathbb{R}\text{ such that }\phi^{\mu}=\pi^{\mu}%
\circ\phi, \label{GCS.1}%
\end{equation}
where $\pi^{\mu}$ are the projection mappings of $\mathbb{R}^{n}$, are the
\emph{coordinate scalar fields} of $(U,\phi).$ For each $p\in U,$ the $n$ real
numbers
\begin{equation}
x^{\mu}=\phi^{\mu}(p) \label{GCS.2}%
\end{equation}
are the so-called \emph{position coordinates} of $p$ with respect to
$(U,\phi).$

The \emph{open subset }$\mathcal{U}^{\prime}$ of the canonical space
$\mathcal{U}_{o}$ defined by
\begin{equation}
\mathcal{U}^{\prime}=\left\{  \lambda^{\mu}b_{\mu}\in\mathcal{U}_{o}\text{
}|\text{ }\lambda^{\mu}\in\phi^{\mu}(U),\text{ for each }\mu=1,\ldots
,n\right\}  \label{GCS.3}%
\end{equation}
will be called a \emph{position vector set} of $U.$ Of course, it is
associated to $(U,\phi).$

There exists an isomorphism between $U$ and $\mathcal{U}^{\prime}$ which is
realized by $U\ni p\mapsto\iota(p)\in\mathcal{U}^{\prime}$ and $\mathcal{U}%
^{\prime}\ni x\mapsto\iota^{-1}(x)\in U$ such that
\begin{align}
\iota(p)  &  =\phi^{\mu}(p)b_{\mu},\label{GCS.4a}\\
\iota^{-1}(x)  &  =\phi^{-1}(b^{1}\cdot x,\ldots,b^{n}\cdot x). \label{GCS.4b}%
\end{align}
As suggested by the above notation , $\iota$ and $\iota^{-1}$ are inverse
mappings of each other.

Indeed, we have that for any $p\in U$
\begin{align*}
\iota^{-1}\circ\iota(p)  &  =\phi^{-1}(b^{1}\cdot\phi^{\mu}(p)b_{\mu}%
,\ldots,b^{n}\cdot\phi^{\mu}(p)b_{\mu})\\
&  =\phi^{-1}(\phi^{1}(p),\ldots,\phi^{n}(p))\\
&  =\phi^{-1}\circ\phi(p)=p,
\end{align*}
i.e., $\iota^{-1}\circ\iota=i_{U}$.

And for any $x\in\mathcal{U}^{\prime}$
\begin{align*}
\iota\circ\iota^{-1}(x)  &  =\phi^{\mu}(\phi^{-1}(b^{1}\cdot x,\ldots
,b^{n}\cdot x))b_{\mu}\\
&  =\pi^{\mu}\circ\phi\circ\phi^{-1}(b^{1}\cdot x,\ldots,b^{n}\cdot x)b_{\mu
}\\
&  =(b^{\mu}\cdot x)b_{\mu}=x,
\end{align*}
i.e., $\iota\circ\iota^{-1}=i_{\mathcal{U}^{\prime}}.$

Then, any \emph{point} $p\in U$ can be also localized by a \emph{vector}
$\iota(p)\in\mathcal{U}^{\prime}.$ We call
\begin{equation}
x=\iota(p),\text{ i.e., }x=x^{\mu}b_{\mu}, \label{GCS.5}%
\end{equation}
the \emph{position vector} of $p$ with respect to $(U,\phi)$.

Let $X$ be a smooth multivector field on $U$. Their position coordinate and
vector representations with respect to $(U,\phi)$ are respectively given by
\begin{equation}
\phi(U)\ni(x^{1},\ldots,x^{n})\mapsto X\circ\phi^{-1}(x^{1},\ldots,x^{n}%
)\in\bigwedge\mathcal{U}_{o} \label{GCS.5a0}%
\end{equation}
and
\begin{equation}
\iota(U)\ni x\mapsto X\circ\iota^{-1}(x)\in\bigwedge\mathcal{U}_{o}.
\label{GCS.5a00}%
\end{equation}

Let us take $a\in\mathcal{U}_{o}.$ For any smooth multivector field $X,$ the
$a$-directional ordinary derivative of the position vector representation of
$X$ with respect to $(U,\phi),$ namely $a\cdot\partial_{x}X\circ\iota^{-1},$
is related to the $x^{\mu}$-partial derivatives of the position coordinate
representation of $X$ with respect to $(U,\phi),$ namely $\dfrac{\partial
}{\partial x^{\mu}}X\circ\phi^{-1},$ by the identity
\begin{equation}
(a\cdot\partial_{x}X\circ\iota^{-1})\circ\iota=a\cdot b^{\mu}(\dfrac{\partial
}{\partial x^{\mu}}X\circ\phi^{-1})\circ\phi, \label{GCS.5a}%
\end{equation}
whose proof is a simple computation.

\subsection{Enter $a\cdot\partial$}

Let $X$ be a smooth multivector field on $U,$ and let us take $a\in
\mathcal{U}_{o}.$ The smooth multivector field on $U,$ namely $a\cdot\partial
X,$ defined as $U\ni p\mapsto a\cdot\partial X(p)\in\bigwedge\mathcal{U}_{o}$
such that
\begin{equation}
a\cdot\partial X(p)=(a\cdot\partial_{x}X\circ\iota^{-1})(x), \label{GCS.5b}%
\end{equation}
where $x=\iota(p),$ is called the $a$-\emph{directional ordinary derivative}
($a$\emph{-DOD}) of $X$ with respect to $(U,\phi).$

In agreement with Eq.(\ref{GCS.5a}), we have a noticeable formula. For any
$a\in\mathcal{U}_{o},$ and for all $X\in\mathcal{M}(U)$
\begin{equation}
a\cdot\partial X(p)=a\cdot\partial_{x}X\circ\iota^{-1}(x)=a\cdot b^{\mu}%
\dfrac{\partial}{\partial x^{\mu}}X\circ\phi^{-1}(x^{1},\ldots,x^{n}),
\label{GCS.5c}%
\end{equation}
where $x=\iota(p)$ and $(x^{1},\ldots,x^{n})=\phi(p).$

Let us take $a\in\mathcal{V}(U).$ The $a$-\emph{DOD }of $X$ with respect to
$(U,\phi),$ which as usual will be also denoted by $a\cdot\partial X,$ is
defined to be $U\ni p\mapsto a\cdot\partial X(p)\in\bigwedge\mathcal{U}_{o}$
such that
\begin{equation}
a\cdot\partial X(p)=a\circ\iota^{-1}(x)\cdot\partial_{x}X\circ\iota^{-1}(x),
\label{GCS.6a}%
\end{equation}
where $x=\iota(p).$ In accordance with Eq.(\ref{GCS.5a}), we might be also
written
\begin{equation}
a\cdot\partial X(p)=a\circ\phi^{-1}(x^{1},\ldots,x^{n})\cdot b^{\mu}%
\frac{\partial}{\partial x^{\mu}}X\circ\phi^{-1}(x^{1},\ldots,x^{n}),
\label{CGS.6b}%
\end{equation}
where $(x^{1},\ldots,x^{n})=\phi(p)$.

We notice that the basic properties which we expected should be valid for
$a\cdot\partial$ are completely analogous to those ones which are satisfied by
$a\cdot\partial_{o}$.

\subsection{Covariant and Contravariant Frame Fields}

The canonical position vector of any $p\in U,$ namely $x_{o}\in\iota_{o}(U),$
can be considered either as function of its position vector with respect
$(U,\phi),$ namely $x\in\iota(U),$ i.e.,
\begin{equation}
\iota(U)\ni x\mapsto x_{o}=\phi_{o}^{\nu}\circ\iota^{-1}(x)b_{\nu}\in\iota
_{o}(U), \label{GCS.6}%
\end{equation}
or as function of its position coordinates with respect to $(U,\phi),$ namely
$(x^{1},\ldots,x^{n})\in\phi(U),$ i.e.,
\begin{equation}
\phi(U)\ni(x^{1},\ldots,x^{n})\mapsto x_{o}=\phi_{o}^{\nu}\circ\phi^{-1}%
(x^{1},\ldots,x^{n})b_{\nu}\in\iota_{o}(U). \label{GCS.7}%
\end{equation}
We emphasize that the first one is the position vector representation of
$\phi_{o}^{\nu}b_{\nu},$ and the second one is the coordinate vector
representation of $\phi_{o}^{\nu}b_{\nu},$ both of them with respect to
$(U,\phi).$

The $n$ smooth vector fields on $U,$ namely $b_{1}\cdot\partial\phi_{o}^{\nu
}b_{\nu},\ldots,b_{n}\cdot\partial\phi_{o}^{\nu}b_{\nu},$ defines a frame
field on $U$ which is called the \emph{covariant frame field} for $(U,\phi).$
It will be denoted by $\{b_{\mu}\cdot\partial x_{o}\},$ i.e.,
\begin{equation}
b_{\mu}\cdot\partial x_{o}\equiv b_{\mu}\cdot\partial\phi_{o}^{\nu}b_{\nu
},\text{ for each }\mu=1,\ldots,n. \label{GCS.8a}%
\end{equation}
In agreement to Eq.(\ref{GCS.5c}), this means that for any $p\in U$
\begin{equation}
b_{\mu}\cdot\partial x_{o}(p)=b_{\mu}\cdot\partial_{x}\phi_{o}^{\nu}\circ
\iota^{-1}(x)b_{\nu}=\frac{\partial}{\partial x^{\mu}}\phi_{o}^{\nu}\circ
\phi^{-1}(x^{1},\ldots,x^{n})b_{\nu}, \label{GCS.8b}%
\end{equation}
where $x=\iota(p)$ and $(x^{1},\ldots,x^{n})=\phi(p).$

The position coordinates of any $p\in U$ with respect to $(U,\phi),$ namely
$x^{\mu}\in\phi^{\mu}(U)\subseteq\mathbb{R},$ can be taken either as functions
of its \emph{canonical position vector} $x_{o}\in\iota_{o}(U)\subseteq
\mathcal{U}_{o}^{\prime},$ i.e.,
\begin{equation}
\iota_{o}(U)\ni x_{o}\mapsto x^{\mu}=\phi^{\mu}\circ\iota_{o}^{-1}(x_{o}%
)\in\phi^{\mu}(U), \label{GCS.9}%
\end{equation}
or as function of its \emph{canonical position coordinates} $(x_{o}^{1}%
,\ldots,x_{o}^{n})\in\phi_{o}(U)\subseteq\phi_{o}(U_{o}),$ i.e.,
\begin{equation}
\phi_{o}(U)\ni(x_{o}^{1},\ldots,x_{o}^{n})\mapsto x^{\mu}=\phi^{\mu}\circ
\phi_{o}^{-1}(x_{o}^{1},\ldots,x_{o}^{n})\in\phi^{\mu}(U). \label{GCS.10}%
\end{equation}
It should be noted that the former is the \emph{canonical position vector
representation} of $\phi^{\mu}$, meanwhile the latter is the \emph{canonical
position coordinate representation} of $\phi^{\mu}.$

The $n$ smooth vector fields on $U,$ namely $\partial_{o}\phi^{1}%
,\ldots,\partial_{o}\phi^{n},$ defines a frame field on $U$ which is called
the \emph{contravariant frame field} for $(U,\phi).$ It is usually denoted by
$\{\partial_{o}x^{\nu}\}$, where $\partial_{o}$ is the Hestenes operator
introduced in Section 3.3. We have, of course, $\partial_{o}x^{\nu}%
\equiv\partial_{o}\phi^{\nu}$, for each $\nu=1,\ldots,n$, which according to
Eq.(\ref{GCS.10}) means that for any $p\in U$
\begin{equation}
\partial_{o}x^{\nu}(p)=\partial_{x_{o}}\phi^{\nu}\circ\iota_{o}^{-1}%
(x_{o})=b^{\mu}\frac{\partial}{\partial x^{\mu}}\phi^{\nu}\circ\phi_{o}%
^{-1}(x_{o}^{1},\ldots,x_{o}^{n}), \label{GCS.11b}%
\end{equation}
where $x_{o}=\iota_{o}(p)$ and $(x_{o}^{1},\ldots,x_{o}^{n})=\phi_{o}(p).$

\subsection{The Relation Between $a\cdot\partial_{o}$ and $a\cdot\partial$}

To find the relation between the operators $a\cdot\partial_{o}$ and
$a\cdot\partial$ we now introduce a \emph{non-singular} smooth $(1,1)$%
-extensor field on $U,$ namely $J_{\phi},$ defined as
\begin{equation}
J_{\phi}(a)=a\cdot\partial x_{o},\text{ i.e., }J_{\phi}(a)=a\cdot\partial
\phi_{o}^{\nu}b_{\nu}, \label{GCS.12a}%
\end{equation}
which will be called the \emph{Jacobian field} for $(U,\phi)$. Of course, for
any $p\in U$
\begin{equation}
\left.  J_{\phi}\right\vert _{(p)}(a)=a\cdot\partial_{x}\phi_{o}^{\nu}%
\circ\iota^{-1}(x)b_{\nu}=a\cdot b^{\mu}\frac{\partial}{\partial x^{\mu}}%
\phi_{o}^{\nu}\circ\phi^{-1}(x^{1},\ldots,x^{n})b_{\nu}. \label{GCS.12b}%
\end{equation}

The inverse of $J_{\phi}$ denoted by $J_{\phi}^{-1}$ is given by
\begin{equation}
J_{\phi}^{-1}(a)=a\cdot\partial_{o}x,\text{ i.e., }J_{\phi}^{-1}%
(a)=a\cdot\partial_{o}\phi^{\nu}b_{\nu}. \label{GCS.13a}%
\end{equation}
Therefore, for any $p\in U$,
\begin{equation}
\left.  J_{\phi}^{-1}\right\vert _{(p)}(a)=a\cdot\partial_{x_{o}}\phi^{\nu
}\circ\iota_{o}^{-1}(x_{o})b_{\nu}=a\cdot b^{\mu}\frac{\partial}{\partial
x_{o}^{\mu}}\phi^{\nu}\circ\phi_{o}^{-1}(x_{o}^{1},\ldots,x_{o}^{n})b_{\nu}.
\label{GCS.13b}%
\end{equation}

We now show that there exists a bijective \emph{vector-valued function of a
vector variable }$\varphi:\iota_{o}(U)\rightarrow\iota(U),$ and that there are
$n$ injective \emph{scalar-valued functions of a vector variable}
$\varphi^{\mu}:\iota_{o}(U)\rightarrow\phi^{\mu}(U)$ such that
\begin{align}
\varphi^{-1}(x)  &  =\phi_{o}^{\nu}\circ\iota^{-1}(x)b_{\nu},\label{GCS.14a}\\
\varphi^{\nu}(x_{o})  &  =\phi^{\nu}\circ\iota_{o}^{-1}(x_{o}),\label{GCS.14b}%
\\
\varphi^{\nu}\circ\varphi^{-1}(x)  &  =b^{\nu}\cdot x. \label{GCS.14c}%
\end{align}
Such $\varphi$ and $\varphi^{\nu}$ are given by
\begin{align}
\varphi &  =\iota\circ\iota_{o}^{-1},\label{GCS.15a}\\
\varphi^{\nu}(x_{o})  &  =b^{\nu}\cdot\varphi(x_{o}). \label{GCS.15b}%
\end{align}

Indeed, by putting Eq.(\ref{PV.2a}) into $\varphi^{-1}$ as obtained from
Eq.(\ref{GCS.15a}), we get Eq.(\ref{GCS.14a}). By putting Eq.(\ref{GCS.4a})
into Eq.(\ref{GCS.15a}), and then the result obtained into Eq.(\ref{GCS.15b}),
we get Eq.(\ref{GCS.14b}). From Eq.(\ref{GCS.14b}) and $\varphi^{-1}$ as
obtained from Eq.(\ref{GCS.15a}), by using Eq.(\ref{GCS.4b}) and
Eq.(\ref{GCS.1}), it follows Eq.(\ref{GCS.14c}).

We have that $\varphi^{\nu}$ and $\varphi^{-1}$ are involved into remarkable
formulas for the covariant and contravariant frame fields of $(U,\phi)$%
\begin{align}
b_{\mu}\cdot\partial x_{o}(p)  &  =b_{\mu}\cdot\partial_{x}\varphi
^{-1}(x),\label{GCS.16a}\\
\partial_{o}x^{\nu}(p)  &  =\partial_{x_{o}}\varphi^{\nu}(x_{o}).
\label{GCS.16b}%
\end{align}
There are also noticeable formulas for the Jacobian field, and its inverse, in
which $\varphi^{-1}$ and $\varphi$ are involved. They are,
\begin{align}
\left.  J_{\phi}\right|  _{(p)}(a)  &  =a\cdot\partial_{x}\varphi
^{-1}(x),\label{GCS.17a}\\
\left.  J_{\phi}^{-1}\right|  _{(p)}(a)  &  =a\cdot\partial_{x_{o}}%
\varphi(x_{o}), \label{GCS.17b}%
\end{align}
where $x=\iota(p)$ and $x_{o}=\iota_{o}(p).$

Eq.(\ref{GCS.16a}) and Eq.(\ref{GCS.16b}) follow immediately from
Eq.(\ref{GCS.8b}) and Eq.(\ref{GCS.11b}) once we use Eq.(\ref{GCS.14a}) and
Eq.(\ref{GCS.14b}), respectively. Eq.(\ref{GCS.17a}) and Eq.(\ref{GCS.17b})
are immediate consequences of Eq.(\ref{GCS.12b}) and Eq.(\ref{GCS.13b}) once
we take into account Eq.(\ref{GCS.14a}) and Eq.(\ref{GCS.14b}), respectively.

We end this subsection presenting some useful properties, which can be easily
proved and which will be used in the sequel papers of this series.

\textbf{i.} $\{b_{\mu}\cdot\partial x_{o}\}$ and $\{\partial_{o}x^{\mu}\}$
define a pair of \textit{reciprocal frame fields} on $U,$ i.e.,
\begin{equation}
b_{\mu}\cdot\partial x_{o}\cdot\partial_{o}x^{\nu}=\delta_{\mu}^{\nu}.
\label{GCS.18}%
\end{equation}

\textbf{iii.} $\{b_{\mu}\cdot\partial x_{o}\}$ and $\{\partial_{o}x^{\mu}\}$
are just $J_{\phi}$-deformations of the fiducial frame field $\{b_{\mu}\},$
i.e.,
\begin{align}
b_{\mu}\cdot\partial x_{o}  &  =J_{\phi}(b_{\mu}),\label{GCS.20a}\\
\partial_{o}x^{\nu}  &  =J_{\phi}^{\ast}(b^{\nu}). \label{GCS.20b}%
\end{align}

The first statement is trivial from the definition of deformation, see
\cite{3}. In order to prove the second statement, we write
\begin{align*}
a\cdot J_{\phi}^{\ast}(b^{\nu})  &  =J_{\phi}^{-1}(a)\cdot b^{\nu}\\
&  =(a\cdot\partial_{o}x)\cdot b^{\nu}=a\cdot\partial_{o}(x\cdot b^{\nu
})=a\cdot\partial_{o}x^{\nu}=a\cdot(\partial_{o}x^{\nu}),
\end{align*}
hence, by non-degeneracy of scalar product, the required result immediately
follows.\bigskip

\noindent\textbf{Proposition 3. }\textit{The relationship between the
canonical} $a$-\emph{DODO} \textit{and the} $a$-\emph{DODO} \noindent
\textit{with respect to} $(U,\phi),$ \textit{namely} $a\cdot\partial_{o}$
\textit{and} $a\cdot\partial,$ \textit{is given by}
\begin{equation}
J_{\phi}(a)\cdot\partial_{o}X=a\cdot\partial X, \label{GCS.21}%
\end{equation}
\noindent\textit{for all }$X\in\mathcal{M}(U)$.

\noindent\textbf{Proof}

Let $X$ be a smooth multivector field on $U.$ Then, by using a link rule for
the $a$ directional derivation of multivector-valued functions of a vector
variable on the obvious identity $(X\circ\iota_{o}^{-1})\circ\varphi
^{-1}(x)=X\circ\iota^{-1}(x),$ we can write
\begin{align*}
a\cdot\partial_{x}(X\circ\iota_{o}^{-1})\circ\varphi^{-1}(x)  &
=a\cdot\partial_{x}X\circ\iota^{-1}(x)\\
(a\cdot\partial_{x}\varphi^{-1}(x)\cdot\partial_{x_{o}}X\circ\iota_{o}%
^{-1})\circ\varphi^{-1}(x)  &  =a\cdot\partial_{x}X\circ\iota^{-1}(x)\\
a\cdot\partial_{x}\varphi^{-1}(x)\cdot\partial_{x_{o}}X\circ\iota_{o}%
^{-1}(x_{o})  &  =a\cdot\partial_{x}X\circ\iota^{-1}(x).
\end{align*}
Hence, using Eq.(\ref{GCS.17a}), and Eq.(\ref{AOD.3a}) and Eq.(\ref{GCS.5c}),
we get
\[
\left.  J_{\phi}\right\vert _{(p)}(a)\cdot\partial_{o}X(p)=a\cdot\partial
X(p).\blacksquare
\]

We end this section by presenting a remarkable and useful property which is
immediate consequence of Eq.(\ref{GCS.21}) and Eq.(\ref{GCS.20a}). It is
\[
(b_{\mu}\cdot\partial x_{o})\cdot\partial_{o}X=b_{\mu}\cdot\partial X,
\]
for all $X\in\mathcal{M}(U).$

\section{Conclusions}

We just presented a theory of multivector and extensor fields living on a
smooth manifold $M$ of arbitrary topology based on the geometric algebra of
multivectors and extensors. Our approach, which does not suffer the problems
of earlier attempts\footnote{Which are indeed restricted to vector manifolds.}
of applying geometric algebra to the differential geometry of manifolds is
based on calculations done with the representatives of those fields through
canonical algebraic structures over the canonical space $\mathcal{U}_{o}$
associated to a local chart $(U_{o},\phi_{o})$ of the maximal atlas of $M$.
Some crucial ingredients for our program, like the $a$-directional ordinary
derivative of multivector and extensor fields, the Lie algebra of smooth
vector fields and the Hestenes derivatives have been introduced and their main
properties explored. We worked also some few applications involving general
coordinate systems besides the one used in the definition of the canonical
space. \medskip

\textbf{Acknowledgments: }V. V. Fern\'{a}ndez and A. M. Moya are very grateful
to Mrs. Rosa I. Fern\'{a}ndez who gave to them material and spiritual support
at the starting time of their research work. This paper could not have been
written without her inestimable help. Authors are also grateful to Drs. E.
Notte-Cuello and E. Capelas de Oliveira for useful discussions.

\end{document}